\documentclass[10pt]{article}
\usepackage{amssymb,amsmath,amsthm}
\usepackage{enumerate}

\newcommand*{\clo}{\overline}
\newcommand*{\Clo}{\overline}
\newcommand*{\reals}{\ensuremath{\mathbb{R}}}
\newcommand*{\nat}{\ensuremath{\mathbb{N}}}

\newcommand*{\ratio}{\mathbb{Q}}

\renewcommand*{\epsilon}{\varepsilon}
\newcommand*{\conti}{\ensuremath{\mathfrak{c}}}
\newtheorem{theorem}{Theorem}

\newtheorem{lemma}[theorem]{Lemma}

\title{A Nonseparably Connected Metric Space as a Dense Connected Graph}
\begin{document}

\author{Micha{\l} Morayne and Micha{\l} Ryszard W\'{o}jcik
\footnote{This material comes from M.~R.~W\'{o}jcik's
PhD thesis supervised by M.~Morayne, see \cite{PhD}.
It was also presented at a conference
(http://atlas-conferences.com/cgi-bin/abstract/caxg-08).}}

\maketitle
\begin{abstract}
We present a connected metric space
that does not contain any
nontrivial separable connected subspace.
Our space is a dense connected graph of a function
from the real line satisfying Cauchy's equation.
\end{abstract}

A space is {\em separably connected} if any two of its points
can be contained in a connected separable subset.
Clearly, this is a generalization of path connectedness.
This notion has emerged recently in two independent lines of research
--- see \cite{Candeal} and section 4.5 in \cite{PhD}.
We define a {\em nonseparably connected} space to be a connected space
whose all separable connected subsets are singletons.
The first example of a nonseparably connected metric space
was given by Pol \cite{POL} in 1975.
There is also a recent example given by Simon \cite{SIMON} in 2001.
In this note we obtain such a space
as a graph of a function from the real line
satisfying Cauchy's equation.
Consequently, our space is homogeneous,
not locally connected at any point,
and all its points are cut points.

Having observed that a function with a dense connected graph
from the real line into a nonseparable metric space
would serve as an example of a nonseparably connected metric space,
it is natural to review the existing constructions of
dense connected graphs to see if they can be made to work in this context.

In \cite{Phillips} Phillips constructed a dense connected graph in the plane.
Then Kulpa in \cite{Kulpa} generalized this construction and revealed that
the assumption $\text{weight}(Y)\leq\text{weight}(X)$ is used to obtain
a dense connected graph of a function from $X$ to $Y$.
So with $X$ the real line, we end up with a separable graph.
However, the method of obtaining a dense connected graph in the plane due to Jones
\cite{Jones} can be adapted to our purposes.

In fact, we have found a way of constructing a function $f\colon X\to Y$
with a connected dense graph for a broad class of products $X\times Y$,
including any normed spaces $X$ and $Y$ of size $\conti$,
in which case such a function may even be required
to satisfy Cauchy's equation $f(x+u)=f(x)+f(u)$.
\\\\
{\bf Notation.} If $E\subset X\times Y$ then $dom(E)$ is the projection of $E$ onto $X$.
We use both $|X|$ and $card(X)$ to denote the cardinality of $X$.
When $card(X)=\kappa$, we often say that $X$ has size $\kappa$ or is of size $\kappa$.
We denote the power set of $X$ as $\mathcal{P}(X)$.
Since we consider the graph of a function $F\colon X\to Y$
as a subset of the topological space $X\times Y$
or even as a topological space in its own right,
we decided to use the capital letter $F$
rather than the usual $f$.
When we write that $F$ is connected, we mean that $F$ is a connected
subset of $X\times Y$. Sometimes authors write $f$ is connected
to mean that $f$ maps connected sets onto connected sets
instead of the unambiguous $f$ is Darboux, but we never do that.
For a topological space $X$,
$X$ is $T_1$ if all singletons are closed,
and $X$ is $T_5$ if for any $A,B\subset X$ with
$\Clo{A}\cap B = A\cap\Clo{B} = \emptyset$, 
there are disjoint open sets
$U_1, U_2 \subset X$ with $A\subset U_1$ and $B\subset U_2$.
Notice that $T_5$ does not imply $T_1$.
\\\\
The following technical lemmas
reveal all the topological and combinatorial details
that are essential for our construction of dense connected graphs.

\begin{lemma}
\label{transfinite_construction}
Let $X,Y$ and $\mathcal{H}\subset\mathcal{P}(X\times Y)$
be arbitrary sets such that
$|dom(K)|\geq|\mathcal{H}|$ for every $K\in\mathcal{H}$.
Then there exists a function $f\colon X\to Y$
which intersects every member of the family $\mathcal{H}$.
\end{lemma}
\begin{proof}
Let the ordinal $\Gamma=card(\mathcal{H})$ be used to well-order the family
$\mathcal{H}=\{K_\alpha\colon \alpha<\Gamma\}$.
We will be constructing the desired function by producing one by one
elements of its graph $(x_\alpha,f(x_\alpha))$ by transfinite induction over $\alpha<\Gamma$.
In the first step, pick some $(x_0,f(x_0))\in K_0$.
Now, given an ordinal $\beta<\Gamma$, notice that the set
$dom(K_\beta)\setminus\{x_\alpha\colon\alpha<\beta\}$
is not empty because of our cardinality constraint.
So we can safely choose some $x_\beta$ and some corresponding $f(x_\beta)$
with $(x_\beta,f(x_\beta))\in K_\beta$.
Should the set $X\setminus\{x_\alpha\colon\alpha<\Gamma\}$ remain nonempty,
fill the graph of our function in an arbitrary way just to extend the domain to the whole $X$.
\end{proof}

It is worth noting that in the following lemma $Y$ is assumed to be separably connected,
which is a nice advertisement for this little known topological property.

\begin{lemma}
\label{connected_graph}
Let $X$ be a connected space of size $\kappa$ whose each nonempty open subset
contains a closed separable set of size $\kappa$.
Let $Y$ be a separably connected $T_1$ space.
Suppose that $X\times Y$ is a $T_5$ space whose separable subsets are hereditarily separable.
Let $\mathcal{H}$ be the family of all closed separable subsets of $X\times Y$
whose projection on $X$ has size $\kappa$.
Let the function $F\colon X\to Y$ intersect every member of $\mathcal{H}$.
Then the graph of $F$ is connected and dense in $X\times Y$.
\end{lemma}
\begin{proof}
To see that $F$ is dense
notice that it intersects every set of the form $E\times\{y\}$
--- where $E$ is a closed separable set of size $\kappa$ ---
which is a set in $\mathcal{H}$ to be found in every nonempty open subset of $X\times Y$.

Suppose that $F$ is disconnected.
Since $X\times Y$ is $T_5$,
there are two disjoint open sets $U,V\subset X\times Y$ such that $F\subset U\cup V$,
$F\cap U\not=\emptyset$, $F\cap V\not=\emptyset$.
Let $A=dom(U)$ and $B=dom(V)$.
These sets, as projections of open sets, are open in $X$.
Since $F$ is a function with domain $X$, for every $x\in X$ there is a $y\in Y$
with $(x,y)\in F\subset U\cup V$, and thus $X=A\cup B$.
Since $F\cap U\not=\emptyset$, $A$ is nonempty. Similarly, $B$ is nonempty.
Since $X$ is connected it follows that $A\cap B$ is nonempty.
Hence there is a point $x_0\in X$ and two points $y_1,y_2\in Y$
with $(x_0,y_1)\in U$ and $(x_0,y_2)\in V$.
Since $U,V$ are open, there is an open set $G$ containing $x_0$
such that $G\times\{y_1\}\subset U$ and $G\times\{y_2\}\subset V$.
By assumption, there is a closed separable set $E\subset G$ of size $\kappa$.
Since $Y$ is separably connected, there is a closed connected separable set
$W\subset Y$ containing both $y_1$ and $y_2$.
Notice that the set $K=(E\times W)\setminus(U\cup V)$
is closed and separable in $X\times Y$.
To show that $E\subset dom(K)$ take any $x\in E$.
Since $(x,y_1)\in U$ and $(x,y_2)\in V$,
the connected set $\{x\}\times W$ intersects both $U$ and $V$,
which are disjoint open sets, hence it cannot be covered by their union.
Thus $(x,y)\in (\{x\}\times W)\setminus(U\cup V)\subset K$ for some $y\in Y$
ensuring that $x\in dom(K)$.
We showed that $dom(K)$ contains the set $E$ of size $\kappa$.
Thus $K\in\mathcal{H}$ and by assumption $F\cap K\not=\emptyset$.
But $F\cap K=\emptyset$.
This contradiction shows that $F$ is connected.
\end{proof}

\begin{lemma}
\label{cardinality_lemma}
Let $Z$ be an arbitrary topological space.
Let $\mathcal{H}$ be the family of all closed separable subsets of $Z$.
Then $|\mathcal{H}|\leq|Z|^{\aleph_0}$.
\end{lemma}
\begin{proof}
Let $A,B$ be two distinct closed separable subsets of $Z$.
Then there are two sequences $a,b\in Z^\nat$ such that $\clo{a(\nat)}=A$ and $\clo{b(\nat)}=B$.
Since $A\not=B$, it follows that $a\not=b$.
This means that to each closed separable set we can assign a unique element of $Z^\nat$,
which completes the argument.
\end{proof}

The following theorem yields the existence of functions
$F\colon X\to Y$ with connected dense graphs
inside a broad class of products $X\times Y$.

\begin{theorem}
\label{kappa}
Let $X$ be a connected space whose each nonempty open subset
contains a closed separable set of size $\conti$.
Let $Y$ be a separably connected $T_1$ space.
Suppose that $X\times Y$ is a $T_5$ space of size $\conti$
whose separable subsets are hereditarily separable.
Then there exists a function $F\colon X\to Y$
whose graph is connected and dense in $X\times Y$.
\end{theorem}
\begin{proof}
Let $\mathcal{H}$ be the family of all closed and separable subsets of $X\times Y$
whose projection on $X$ has  size $\conti$.
By Lemma \ref{cardinality_lemma},
$|\mathcal{H}|\leq\conti^{\aleph_0}=\conti$.
Thus $X,Y,\mathcal{H}$ satisfy the assumptions of
Lemma \ref{transfinite_construction}.
So there exists a function $F\colon X\to Y$
which intersects every member of $\mathcal{H}$.
By Lemma \ref{connected_graph},
the graph of $F$ is connected and dense in $X\times Y$.
\end{proof}

In the next theorem we show that if $X,Y$ are normed spaces,
$F$ may satisfy Cauchy's equation.

\begin{theorem}\label{cauchy}
Let $X$, $Y$ be metrizable topological vector spaces
of size $\conti$.
Then there exists a function $F\colon X\to Y$
satisfying $F(x+u)=F(x)+F(u)$ for all $x,u\in X$
whose graph is connected and dense in $X\times Y$.
\end{theorem}
\begin{proof}
Let $\mathcal{H}$ be the family of all closed and separable subsets of $X\times Y$
whose projection on $X$ has  size $\conti$.
By Lemma \ref{cardinality_lemma}, $|\mathcal{H}|\leq\conti^{\aleph_0}=\conti$.
If $E\subset X$ then let $\text{span}_{\ratio}(E)=
\{\Sigma_{i=1}^{n}a_i x_i\colon a_i\in\ratio,\ x_i\in E,\ i\in[1,n]\cap\nat,\ n\in\nat\}$.
Notice that $|E|<\conti\implies|\text{span}_{\ratio}(E)|<\conti$.
Let us write $\mathcal{H}=\{K_\alpha\colon \alpha<\conti\}$.
By transfinite induction over $\alpha<\conti$ we will construct a linearly independent
subset $B_0$ of vectors in $X$ over the field of rational numbers
and a function $F\colon B_0\to Y$ intersecting every member of $\mathcal{H}$.
In the step zero, we choose a nonzero vector $x_0$ such that $(x_0,F(x_0))\in K_0$.
In the $\beta$th step, we choose an $x_\beta\in dom(K_\beta)\setminus\text{span}_\ratio(\{x_\alpha\colon\alpha<\beta\})$
and some $F(x_\beta)$ such that $(x_\beta,F(x_\beta))\in K_\beta$.
Let $B_0=\{x_\alpha\colon\alpha<\conti\}$.
Let us extend $F$ to the whole $X$ in the following way.
Let $B$ be a Hamel basis for $X$ over the field of rational numbers containing the linearly independent set $B_0$.
Put $F(x)=0$ for all $x\in B\setminus B_0$ and put
$F(\Sigma_{i=1}^{n}a_i x_i)=\Sigma_{i=1}^{n}a_i F(x_i)$ for all
$a_i\in\ratio,\ x_i\in B,\ i\in[1,n]\cap\nat,\ n\in\nat$.
Now, the function $F$ satisfies $F(x+u)=F(x)+F(u)$ and by Lemma \ref{connected_graph}
its graph is connected and dense in $X\times Y$.
\end{proof}

We are now ready to construct our connected metric space
whose all nontrivial connected subsets are nonseparable
%whose all connected separable subsets are singletons
(= nonseparably connected).

\begin{theorem}
There exists a nonseparably connected metric space $\mathrm{M}$ of size $\conti$
with the following properties:
\begin{enumerate}[(1)]
\item $\mathrm{M}\setminus\{p\}$ is disconnected for every $p\in\mathrm{M}$,
\item $\mathrm{M}$ is not locally connected at any point,
\item $\mathrm{M}$ is a topological group and thus a homogeneous space.
\end{enumerate}
\end{theorem}
\begin{proof}
Let $Y=l^{\infty}$ be the vector space of all bounded sequences of real numbers,
with the supremum norm.
Clearly $Y$ is nonseparable and of size $\conti$.
By Theorem \ref{cauchy}, we obtain a function $\mathrm{M}\colon\reals\to Y$
satisfying Cauchy's equation
whose graph is connected and dense in $\reals\times Y$.
Let $E$ be a connected subset of $\mathrm{M}$ containing two distinct points,
say $a,b\in dom(E)$ with $a<b$.
Then the set $E_0=\mathrm{M}\cap([a,b]\times Y)$
is contained in $E$ because otherwise $E$ would be disconnected,
separated by $(-\infty,c)\times Y$ and $(c,\infty)\times Y$
for some $c\in(a,b)$.
Since $E_0$ is dense in $[a,b]\times Y$,
which is nonseparable, $E$ is nonseparable, too.

That every point of $\mathrm{M}$ is a cut point --- 
$\mathrm{M}\setminus\{p\}$ is disconnected for every $p\in\mathrm{M}$ ---
follows naturally from the fact that the domain of our function is the real line.
Being discontinuous everywhere,
$\mathrm{M}$ cannot be locally connected at any point
because of Theorem 4 in \cite{PhD}.
Finally, $\mathrm{M}$ is easily seen to be a topological group with the addition operation
given by $(x,\mathrm{M}(x))+(u,\mathrm{M}(u))=(x+u,\mathrm{M}(x)+\mathrm{M}(u))=
(x+u,\mathrm{M}(x+u))$.
\end{proof}

{\bf Problem 1.}
Now we know of three examples of nonseparably connected metric spaces.
Since they are all constructed with the use of the axiom of choice,
we would like to ask whether such a space exists in ZF.

{\bf Problem 2.}
There are a number of connected metric spaces which fail to be
separably connected but contain many nontrivial connected separable subsets.
Two examples are already published by Aron and Maestre in \cite{ARON}
and by Simon in \cite{SIMON},
and two more are not yet published \cite{Taras} and \cite{DF}.
The space given in \cite{Taras} is complete.
We would like to ask whether there exists a complete nonseparably connected metric space
(or at least a Borel subset of a complete metric space).

We wish to thank professors Aleksander B{\l}aszczyk,
Pawe{\l} Krupski, and Wies{\l}aw Kulpa for valuable comments.

\end{document}